\documentclass [a4paper,12pt, reqno]{amsart}
\usepackage{amsfonts}
\usepackage{amssymb}
\usepackage{latexsym}
\usepackage[mathscr]{eucal}
\usepackage[T1]{fontenc}
\usepackage[utf8]{inputenc}
\usepackage{latexsym}
\usepackage{amsmath}
\usepackage{amsthm}

\usepackage{tikz}
\usetikzlibrary{datavisualization}
\usetikzlibrary{datavisualization.formats.functions}

\usepackage{float}

\numberwithin{equation}{section}

\begin{document}

\title[Inequalities of the Edmundson-Lah-Ribari\v c type]{Some inequalities of the Edmundson-Lah-Ribari\v c type for 3-convex functions with applications}
\author[R. Miki\' c]{Rozarija Miki\' c}

\address{%
Faculty of Textile Technology, University of Zagreb\\
Prilaz baruna Filipovi\' ca 28a\\
10 000 Zagreb\\
Croatia}

\email{rozarija.jaksic@ttf.hr}


\author[\dj . Pe\v cari\' c]{\dj ilda Pe\v cari\' c}
\address{Catholic University of Croatia\\
Ilica 242\\
10 000 Zagreb\\
Croatia}
\email{gildapeca@gmail.com}

\author[J. Pe\v cari\' c]{Josip Pe\v cari\' c}
\address{Faculty of Textile Technology, University of Zagreb\\
Prilaz baruna Filipovi\' ca 28a\\
10 000 Zagreb\\
Croatia}
\email{pecaric@hazu.hr}

\subjclass[2010]{26A16, 60E05, 60E15}

\keywords{Jensen inequality, Edmundson-Lah-Ribari\v c inequality, 3-convex functions, $f$-divergence, Zipf-Mandelbrot law, exponential convexity, Stolarsky-type means}

\maketitle

\begin{abstract}
In this paper we derive some Edmundson-Lah-Ribari\v c type inequalities for positive linear functionals and 3-convex functions. Main results are applied to the generalized $f$-divergence functional. Examples with Zipf–Mandelbrot law are used to illustrate the results. In addition, obtained results are utilized in constructing some families of exponentially convex functions and Stolarsky-type means.
\end{abstract}

\section{Introduction}

\newtheorem{tm}{Theorem}[section]
\newtheorem{kor}{Corollary}[section]
\newtheorem{rem}{Remark}[section]
\newtheorem{primjer}{Example}[section]
\newtheorem{df}{Definition}[section]

The importance of Jensen's inequality for convex functions is its applicability in various
branches of mathematics, especially in mathematical analysis and
statistics. In this paper we refer to a general form of the Jensen inequality for positive linear functionals. In order to present our result, we first need to introduce the appropriate setting.

Let $E$ be a non-empty set and let $L$ be a vector space of real-valued functions $f \colon E \to \mathbb{R}$
having the properties:
\begin{itemize}
\item[(L1)]
 $f,g \in L \Rightarrow (af+bg)\in L$ for all $a,b \in\mathbb{R}$;
\item[(L2)]
 $\boldsymbol{1} \in L$, i.e., if $f(t)=1$ for every $t \in E$, then $f \in L$. 
\end{itemize}

We also consider positive linear functionals $A \colon L \to \mathbb{R}$. That is, we assume that:
\begin{itemize}
\item[(A1)]
$A(af+bg)=aA(f)+bA(g)$ for $f,g \in L$ and $a,b \in \mathbb{R}$;
\item[(A2)]
$f \in L$, $f(t) \ge 0$ for every $t \in E \Rightarrow A(f) \ge 0$ ($A$ is positive). 
\end{itemize}

Since it was proved, the famous Jensen inequality and its converses have been extensively studied by many authors and have been generalized in numerous directions. Jessen \cite{jessen} gave the following generalization of Jensen's inequality for convex functions (see also \cite[p.47]{crvena}):
\begin{tm} $\mathrm{(}$\cite{jessen}$\mathrm{)}$
Let $L$ satisfy properties (L1) and (L2) on a nonempty set $E$, and assume that $\phi$ is a continous convex function 
on an interval $I \subset \mathbb{R}$. If $A$ is a positive linear functional with $A(1)=1$, then for all $f \in L$ such
that $\phi (f) \in L$ we have $A(f) \in I$ and
\begin{equation}\label{jensen}
\phi (A(f)) \le A(\phi (f)).
\end{equation} 
\end{tm}

The following result is one of the most famous converses of the Jensen inequality known as the Edmundson-Lah-Ribari\v c inequality, and it was proved in \cite{beesack} by Beesack and Pe\v{c}ari\'{c} (see also \cite[p.98]{crvena}):
\begin{tm} $\mathrm{(}$\cite{beesack}$\mathrm{)}$
Let $\phi$ be convex on the interval $I=[m,M]$ such that $-\infty <m<M<\infty$. Let $L$ satisfy conditions (L1) and (L2) on $E$ and let $A$ be any
 positive linear functional on $L$ with $A(1)=1$. Then for every $f \in L$ such that $\phi (f) \in L$ (so that
$m \le f(t) \le M$ for all $t \in E$), we have
\begin{equation} \label{1}
A(\phi (f)) \le \frac{M-A(f)}{M-m}\phi (m)+\frac{A(f)-m}{M-m} \phi (M).
\end{equation} 
\end{tm}

For some recent results on the converses of the Jensen inequality, the reader is referred to \cite{hor}, \cite{hor1}, \cite{ivelic}, \cite{drugi}, \cite{treci}, \cite{mon6} and \cite{moh}.

Unlike the results from the above mentioned papers, which require convexity of the involved functions, the main objective of this paper is to derive a class of inequalities of the Jensen and Edmundson-Lah-Ribari\v c type that hold for 3-convex functions, which will be an extension of the results from \cite{prvi3conv}.

Definition of the $n$-convex function is characterized by $n$-th order divided difference.
The $n$-th order divided difference of a function $f\colon [a, b] \to \mathbb{R}$ at mutually distinct points $t_0, t_1, ... , t_n\in  [a, b]$ is defined recursively by
\begin{align*}
[t_i]f&=f(t_i), \ \ i=0,...,n,\\
[t_0,...,t_n]f&=\dfrac{[t_1,..., t_n]f-[t_0,..., t_{n-1}]f}{t_n-t_0}.
\end{align*}
The value $[t_0, . . . , t_n]f$ is independent of the order of the points $t_0, . . . , t_n$.

A function $f\colon [a, b] \to \mathbb{R}$ is said to be $n$-convex ($n \ge 0$)
if and only if for all choices of $(n+1)$ distinct points $t_0, t_1, ... , t_n\in  [a, b]$, we have $[t_0, . . . , t_n]f\ge 0$.

In this paper we study 3-convex functions, which are characterized by the third order divided differences. Definition of the third order divided differences can be extended to include the cases in which some or all the points coincide (see e.g. \cite[p.14]{crvena}):

\begin{itemize}
    \item If the function $f$ is differentiable on $[a,b]$ and $t,t_0,t_1\in [a,b]$ are mutually different points, then
    \begin{align}
        [t,t,t_0,t_1]f&=\dfrac{f'(t)}{(t-t_0)(t-t_1)}+\dfrac{f(t)(t_0+t_1-2t)}{(t-t_0)^2(t-t_1)^2}\nonumber\\
        & \ \ \ \ +\dfrac{f(t_0)}{(t_0-t)^2(t_0-t_1)}+\dfrac{f(t_1)}{(t_1-t)^2(t_1-t_0)}.\label{3conv1}
    \end{align}
    \item If the function $f$ is differentiable on $[a,b]$ and $t,t_0\in [a,b]$ are mutually different points, then
    \begin{align}
        [t,t,t_0,t_0]f&=\dfrac{1}{(t_0-t)^3}\left[(t_0-t)(f'(t_0)+f'(t))+2(f(t)-f(t_0))\right].\label{3conv2}
    \end{align}
    \item If the function $f$ is twice differentiable on $[a,b]$ and $t,t_0\in [a,b]$ are mutually different points, then
    \begin{align}
        [t,t,t,t_0]f&=\dfrac{1}{(t_0-t)^3}\left[f(t_0)-\sum_{k=0}^2\dfrac{f^{(k)}(t)}{k!}(t_0-t)^k\right].\label{3conv3}
    \end{align}
    \item If the function $f$ is three times differentiable on $[a,b]$ and $t\in [a,b]$, then
    \begin{align}
        [t,t,t,t]f&=\dfrac{f'''(t)}{3!}.\label{3conv4}
    \end{align}
\end{itemize}

We can extend the definition of 3-convex functions by including the cases in which some or all of the points coincide. This is given in the following theorem which can be easily proven by using the mean value theorem for divided differences (see e.g. \cite{issa}).

\begin{tm}\label{tm_uv}
Let a function $f$ be defined on an interval $I\subseteq \mathbb{R}$. The following equivalences hold.
\begin{itemize}
    \item[(i)] If $f\in \mathcal{C}(I)$, then $f$ is 3-convex if and only if $[t,t,t_0,t_1]f\ge 0$ for all mutually different points $t,t_0,t_1\in I$.
    \item[(ii)] If $f\in \mathcal{C}(I)$, then $f$ is 3-convex if and only if $[t,t,t_0,t_0]f\ge 0$ for all mutually different points $t,t_0\in I$.
    \item[(iii)] If $f\in \mathcal{C}^2(I)$, then $f$ is 3-convex if and only if $[t,t,t,t_0]f\ge 0$ for all mutually different points $t,t_0\in I$.
    \item[(iv)] If $f\in \mathcal{C}^3(I)$, then $f$ is 3-convex if and only if $[t,t,t,t]f\ge 0$ for every $t\in I$.
\end{itemize}
\end{tm}


\section{Results}\label{sec2}

Throughout this paper, whenever mentioning the interval $[m,M]$, we assume that $-\infty<m<M<\infty$ holds.

\begin{tm}\label{tm1}
Let $L$ satisfy conditions (L1) and (L2) on a non-empty set $E$ and let $A$ be any  positive linear functional on $L$ with $A(\boldsymbol{1})=1$. Let $\phi$ be a 3-convex function defined on an interval of real numbers $I$ whose interior contains the interval $[m,M]$.
Then 
 \begin{align}
&\dfrac{A\left[(M\boldsymbol{1}-f)(f-m\boldsymbol{1})\right]}{M-m}\left(\dfrac{\phi(M)-\phi(m)}{M-m}-\phi_+'(m)\right)\nonumber \\
   \le   &\dfrac{M-A(f)}{M-m}\phi (m)+\dfrac{A(f)-m}{M-m}\phi (M)-A\left(\phi (f)\right)\label{tm1rez1}\\
  \le  &\dfrac{A\left[(M\boldsymbol{1}-f)(f-m\boldsymbol{1})\right]}{M-m}\left(\phi_-'(M)-\dfrac{\phi(M)-\phi(m)}{M-m}\right)\nonumber 
 \end{align}
 holds for any $f \in L$ such that $ \phi \circ f \in L$ and $m\le f(t)\le M$ for $t\in E$. If the function $-\phi$ is 3-convex, then the inequalities are reversed. 
\end{tm}

\emph{Proof:} The function $\phi$ is 3-convex, so from Theorem \ref{tm_uv} (i) we have that $[t,t,t_0,t_1]\phi \ge 0$ for all mutually different points $t,t_0,t_1\in I$. When we take $t=m$, $t_0=x$ and $t_1=M$ in (\ref{3conv1}), we obtain that
 \begin{align*}
     0&\le \dfrac{\phi_+'(m)}{(m-x)(m-M)}+\dfrac{\phi(m)(x+M-2m)}{(m-x)^2(m-M)^2}\nonumber\\
        & \ \ \ \         +\dfrac{\phi(x)}{(x-m)^2(x-M)}+\dfrac{\phi(M)}{(M-m)^2(M-x)}
 \end{align*}
 holds for every $x\in \langle m,M\rangle$.
After multiplying by $(x-m)^2(x-M)$ and rearranging, the upper relation becomes
\begin{align}
      &\dfrac{(M-x)(x-m)}{M-m}\left(\dfrac{\phi(M)-\phi(m)}{M-m}-\phi'_+(m)\right) \nonumber\\
      \le   &\dfrac{M-x}{M-m}\phi (m)+\dfrac{x-m}{M-m}\phi (M)-\phi (x).\label{tm1pom1}
\end{align}
In a similar manner, when we put $t=M$, $t_0=x$ and $t_1=m$ in (\ref{3conv1}), after arranging the relation thus obtained, we get that
\begin{align}
  &\dfrac{M-x}{M-m}\phi (m)+\dfrac{x-m}{M-m}\phi (M)-\phi (x)\nonumber\\
  \le  &\dfrac{(M-x)(x-m)}{M-m}\left(\phi'_-(M)-\dfrac{\phi(M)-\phi(m)}{M-m}\right) \label{tm1pom1a}
\end{align}
holds for every $x\in \langle m,M\rangle$. Now, we see that (\ref{tm1pom1}) and (\ref{tm1pom1a}) give the following sequence of inequalities:
\begin{align}
      &\dfrac{(M-x)(x-m)}{M-m}\left(\dfrac{\phi(M)-\phi(m)}{M-m}-\phi'_+(m)\right) \nonumber\\
\le  &\dfrac{M-x}{M-m}\phi (m)+\dfrac{x-m}{M-m}\phi (M)-\phi (x)\nonumber\\
  \le  &\dfrac{(M-x)(x-m)}{M-m}\left(\phi'_-(M)-\dfrac{\phi(M)-\phi(m)}{M-m}\right). \label{tm1pom2}
\end{align}
Since the function $f\in L$ satisfies the bounds $m\le f(t)\le M$, we can replace $x$ with $f(t)$ in (\ref{tm1pom2}), and get
\begin{align*}
      &\dfrac{(M-f(t))(f(t)-m)}{M-m}\left(\dfrac{\phi(M)-\phi(m)}{M-m}-\phi'_+(m)\right) \nonumber\\
\le  &\dfrac{M-f(t)}{M-m}\phi (m)+\dfrac{f(t)-m}{M-m}\phi (M)-\phi (f(t))\nonumber\\
  \le  &\dfrac{(M-f(t))(f(t)-m)}{M-m}\left(\phi'_-(M)-\dfrac{\phi(M)-\phi(m)}{M-m}\right). 
  \end{align*}
  The inequalities (\ref{tm1rez1}) follow after applying linear functional $A$ to the previous relation taking into account linearity of the functional $A$ and condition $A(\boldsymbol{1})=1$.
\begin{flushright}
 $\Box$
\end{flushright}

\begin{rem}\label{rem2}
The result from Theorem \ref{tm1} is already proven in the paper \cite{prvi3conv}, but in this paper we have provided a shorter and more elegant proof.
\end{rem}

\begin{tm}\label{tm2}
Let $L$ satisfy conditions (L1) and (L2) on a non-empty set $E$ and let $A$ be any  positive linear functional on $L$ with $A(\boldsymbol{1})=1$. Let $\phi$ be a 3-convex function defined on an interval of real numbers $I$ whose interior contains the interval $[m,M]$ and differentiable on $\langle m,M\rangle$. Then 
 \begin{align}
 &(A(f)-m)\left[\dfrac{\phi(M)-\phi(m)}{M-m}-\dfrac{\phi_+'(m)}{2}\right]-\dfrac{1}{2}A[(f-m\boldsymbol{1})\phi'(f)] \nonumber\\
\le  &\dfrac{M-A(f)}{M-m}\phi (m)+\dfrac{A(f)-m}{M-m}\phi (M)-A(\phi (f))\label{tm2rez1}\\
  \le  & \dfrac{1}{2}A[(M\boldsymbol{1}-f)\phi'(f)]- (M-A(f))\left[\dfrac{\phi(M)-\phi(m)}{M-m}-\dfrac{\phi_-'(M)}{2}\right] \nonumber
 \end{align}
 holds for any $f \in L$ such that $ \phi \circ f \in L$ and $m\le f(t)\le M$ for $t\in E$. If the function $-\phi$ is 3-convex, then the inequalities are reversed. 
\end{tm}

\emph{Proof:} Let $\phi$ be a 3-convex function. From Theorem \ref{tm_uv} (ii) we have that $[t,t,t_0,t_0]\phi \ge 0$ for all mutually different points $t,t_0\in I$. When we take $t=m$ and $t_0=x$ in (\ref{3conv2}), we obtain that
 \begin{align*}
     0&\le \dfrac{1}{(x-m)^3}\left[(x-m)(\phi'(x)+\phi_+'(m))+2(\phi(m)-\phi(x))\right]
 \end{align*}
 holds for every $x\in \langle m,M\rangle$.
After multiplying by $(x-m)^3$ and rearranging, the relation from above becomes
\begin{align}
      &(x-m)\left[\dfrac{\phi(M)-\phi(m)}{M-m}-\dfrac{1}{2}\left(\phi'(x)+\phi_+'(m)\right)\right] \nonumber\\
      \le   &\dfrac{M-x}{M-m}\phi (m)+\dfrac{x-m}{M-m}\phi (M)-\phi (x).\label{tm2pom1}
\end{align}
Similarly, when we put $t=M$ and $t_0=x$ in (\ref{3conv2}) and rearrange the obtained relation, we get that
\begin{align}
&\dfrac{M-x}{M-m}\phi (m)+\dfrac{x-m}{M-m}\phi (M)-\phi (x) \nonumber\\
  \le   &(M-x)\left[\dfrac{1}{2}\left(\phi'(x)+\phi_-'(M)\right)-\dfrac{\phi(M)-\phi(m)}{M-m}\right]  \label{tm2pom1a}
\end{align}
holds for every $x\in \langle m,M\rangle$. Now, we see that (\ref{tm2pom1}) and (\ref{tm2pom1a}) together give the following sequence of inequalities:
\begin{align}
      &(x-m)\left[\dfrac{\phi(M)-\phi(m)}{M-m}-\dfrac{1}{2}\left(\phi'(x)+\phi_+'(m)\right)\right] \nonumber\\
\le  &\dfrac{M-x}{M-m}\phi (m)+\dfrac{x-m}{M-m}\phi (M)-\phi (x)\nonumber\\
  \le  &(M-x)\left[\dfrac{1}{2}\left(\phi'(x)+\phi_-'(M)\right)-\dfrac{\phi(M)-\phi(m)}{M-m}\right]  .\label{tm2pom2}
\end{align}
Since the function $f\in L$ satisfies the bounds $m\le f(t)\le M$, we can replace $x$ with $f(t)$ in (\ref{tm2pom2}), and get
\begin{align*}
   &(f(t)-m)\left[\dfrac{\phi(M)-\phi(m)}{M-m}-\dfrac{1}{2}\left(\phi'(f(t))+\phi'_+(m)\right)\right] \nonumber\\
\le  &\dfrac{M-f(t)}{M-m}\phi (m)+\dfrac{f(t)-m}{M-m}\phi (M)-\phi (f(t))\nonumber\\
  \le  &(M-f(t))\left[\dfrac{1}{2}\left(\phi'(f(t))+\phi'_-(M)\right)-\dfrac{\phi(M)-\phi(m)}{M-m}\right]  . 
  \end{align*}
  The inequalities (\ref{tm2rez1}) follow after applying linear functional $A$ to the previous relation taking into account linearity of the functional $A$ and condition $A(\boldsymbol{1})=1$.
\begin{flushright}
 $\Box$
\end{flushright}

\begin{rem}\label{rem1}
If it exists, the first derivative $\phi'$ of a 3-convex function $\phi$ is a convex function. It is known that convex functions are continuous on every open interval, and their one-sided derivatives exist and are finite.
\end{rem}

\begin{tm}\label{tm3}
Let $L$ satisfy conditions (L1) and (L2) on a non-empty set $E$ and let $A$ be any  positive linear functional on $L$ with $A(\boldsymbol{1})=1$. Let $\phi$ be a 3-convex function defined on an interval of real numbers $I$ whose interior contains the interval $[m,M]$ and differentiable on $\langle m,M\rangle$. Then 
 \begin{align}
       &(M-A(f))\left[\phi_-'(M)-\dfrac{\phi(M)-\phi(m)}{M-m}\right]-\dfrac{\phi_-''(M)}{2}A[(M\boldsymbol{1}-f)^2] \nonumber\\
\le  &\dfrac{M-A(f)}{M-m}\phi (m)+\dfrac{A(f)-m}{M-m}\phi (M)-A(\phi (f))\label{tm3rez1}\\
  \le  &(A(f)-m)\left[\dfrac{\phi(M)-\phi(m)}{M-m}-\phi_+'(m)\right]-\dfrac{\phi_+''(m)}{2}A[(f-m\boldsymbol{1})^2] \nonumber
 \end{align}
 holds for any $f \in L$ such that $ \phi \circ f \in L$ and $m\le f(t)\le M$ for $t\in E$. If the function $-\phi$ is 3-convex, then the inequalities are reversed. 
\end{tm}

\emph{Proof:} The function $\phi$ is 3-convex on $[m,M]$ and twice differentiable, so from Theorem \ref{tm_uv} (iii) we have that $[t,t,t_0,t_1]\phi \ge 0$ for all mutually different points $t,t,t_0\in [m,M]$. When we take $t=m$ and $t_0=x$ in (\ref{3conv3}), we obtain that
 \begin{align*}
     0&\le \dfrac{1}{(x-m)^3}\left[\phi(x)-\sum_{k=0}^2\dfrac{\phi_+^{(k)}(m)}{k!}(x-m)^k\right]
 \end{align*}
 holds for every $x\in \langle m,M\rangle$.
After multiplying by $(x-m)^3$ and rearranging, the upper relation becomes
\begin{align}
&\dfrac{M-x}{M-m}\phi (m)+\dfrac{x-m}{M-m}\phi (M)-\phi (x)\nonumber\\
    \le  &(x-m)\left[\dfrac{\phi(M)-\phi(m)}{M-m}-\phi_+'(m)-\dfrac{\phi_+''(m)}{2}(x-m) \right].\label{tm3pom1}
\end{align}
In a similar manner, when we put $t=M$ and $t_0=x$ in (\ref{3conv3}), after rearranging the relation thus obtained, we get that
\begin{align}
   &(M-x)\left[\phi_-'(M)-\dfrac{\phi(M)-\phi(m)}{M-m}-\dfrac{\phi_-''(M)}{2}(M-x) \right]\nonumber\\
    \le   &\dfrac{M-x}{M-m}\phi (m)+\dfrac{x-m}{M-m}\phi (M)-\phi (x) \label{tm3pom1a}
\end{align}
holds for every $x\in \langle m,M\rangle$. Now, we see that (\ref{tm3pom1}) and (\ref{tm3pom1a}) give the following sequence of inequalities:
\begin{align}
      &(M-x)\left[\phi_-'(M)-\dfrac{\phi(M)-\phi(m)}{M-m}-\dfrac{\phi_-''(M)}{2}(M-x) \right]\nonumber\\
\le  &\dfrac{M-x}{M-m}\phi (m)+\dfrac{x-m}{M-m}\phi (M)-\phi (x)\nonumber\\
  \le  &(x-m)\left[\dfrac{\phi(M)-\phi(m)}{M-m}-\phi_+'(m)-\dfrac{\phi_+''(m)}{2}(x-m) \right]. \label{tm3pom2}
\end{align}
Since the function $f\in L$ satisfies the bounds $m\le f(t)\le M$, we can replace $x$ with $f(t)$ in (\ref{tm3pom2}), and get
\begin{align*}
       &(M-f(t))\left[\phi_-'(M)-\dfrac{\phi(M)-\phi(m)}{M-m}-\dfrac{\phi_-''(M)}{2}(M-f(t)) \right]\nonumber\\
\le  &\dfrac{M-f(t)}{M-m}\phi (m)+\dfrac{f(t)-m}{M-m}\phi (M)-\phi (f(t))\nonumber\\
  \le  &(f(t)-m)\left[\dfrac{\phi(M)-\phi(m)}{M-m}-\phi_+'(m)-\dfrac{\phi_+''(m)}{2}(f(t)-m) \right].
  \end{align*}
  The inequalities (\ref{tm3rez1}) follow after applying linear functional $A$ to the previous relation taking into account linearity of the functional $A$ and condition $A(\boldsymbol{1})=1$.
\begin{flushright}
 $\Box$
\end{flushright}

\begin{rem}
Theorems \ref{tm2} and \ref{tm3} can be utilized for obtaining following Jensen-type inequalities for 3-convex functions. 
\begin{itemize}
    \item[(i)] When we put $x=A(f)$ in scalar inequalities (\ref{tm2pom2}) and then subtract the inequalities from Theorem \ref{tm2}, we get
\begin{align*}
      &(A(f)-m)\left[\dfrac{\phi(M)-\phi(m)}{M-m}-\dfrac{1}{2}\left(\phi'(A(f))+\phi_+'(m)\right)\right] \\
     & -\dfrac{1}{2}A[(M\boldsymbol{1}-f)\phi'(f)]- (M-A(f))\left[\dfrac{\phi(M)-\phi(m)}{M-m}+\dfrac{\phi_-'(M)}{2}\right]\\
     \le & A(\phi (f))-\phi (A(f))\le 
     (M-A(f))\left[\dfrac{1}{2}\left(\phi'(A(f))+\phi_-'(M)\right)-\dfrac{\phi(M)-\phi(m)}{M-m}\right] \\
      & -(A(f)-m)\left[\dfrac{\phi(M)-\phi(m)}{M-m}-\dfrac{\phi_+'(m)}{2}\right]+\dfrac{1}{2}A[(f-m\boldsymbol{1})\phi'(f)].
\end{align*}
    \item[(ii)] When we put $x=A(f)$ in scalar inequalities (\ref{tm3pom2}) and then subtract the inequalities from Theorem \ref{tm3}, we get
\begin{align*}
    &(M-A(f))\left[\phi_-'(M)-\dfrac{\phi(M)-\phi(m)}{M-m}-\dfrac{\phi_-''(M)}{2}(M-A(f)) \right]\\
    &-(A(f)-m)\left[\dfrac{\phi(M)-\phi(m)}{M-m}-\phi_+'(m)\right]+\dfrac{\phi_+''(m)}{2}A[(f-m\boldsymbol{1})^2]\\
    \le & A(\phi (f))-\phi (A(f))\le
    (A(f)-m)\left[\dfrac{\phi(M)-\phi(m)}{M-m}-\phi_+'(m)-\dfrac{\phi_+''(m)}{2}(A(f)-m) \right]\\
    & -(M-A(f))\left[\phi_-'(M)-\dfrac{\phi(M)-\phi(m)}{M-m}\right]+\dfrac{\phi_-''(M)}{2}A[(M\boldsymbol{1}-f)^2].
\end{align*}
\end{itemize}
\end{rem}

\begin{rem}
Theorem \ref{tm1} can likewise be used to obtain Jensen-type inequalities for 3-convex functions, and that result is already given in \cite{prvi3conv}.
\end{rem}

\section{Applications to Csisz\' ar divergence and Zipf-Mandelbrot law}

Let us denote the set of all probability distributions by $\mathbb{P}$, that is we say $\boldsymbol{p}=(p_1,...,p_n)\in \mathbb{P}$ if $p_i\in [0,1]$ for $i=1,...,n$ and $\sum_{i=1}^np_i=1$. 

Numerous theoretic divergence measures between two
probability distributions 
have been introduced and comprehensively studied. Their applications can be found in the analysis of contingency tables \cite{10}, in approximation of probability distributions \cite{6}, \cite{21}, in signal processing \cite{14}, and in pattern recognition \cite{3}, \cite{4}.

Csisz\' ar \cite{csiszar1}-\cite{csiszar2} introduced the $f-$divergence functional as
\begin{equation}\label{fdiv}
D_f(\boldsymbol{p},\boldsymbol{q})=\sum_{i=1}^nq_if\left(\frac{p_i}{q_i}\right),
\end{equation}
where $f\colon [0,+\infty \rangle$ is a convex function, and it represent a "distance function" on the set of probability distributions $\mathbb{P}$.

A great number of theoretic divergences are special cases of Csisz\' ar $f$-divergence for different choices of the function $f$.

As in Csisz\' ar \cite{csiszar2}, we interpret undefined expressions by
$$f(0)=\lim_{t \to 0^+}f(t), \ \ 0\cdot f\left(\dfrac{0}{0}\right)=0,$$
$$0\cdot f\left(\dfrac{a}{0}\right)=\lim_{\epsilon \to 0^+}f\left(\dfrac{a}{\epsilon}\right)=a\cdot \lim_{t \to \infty}\dfrac{f(t)}{t}.$$

In this section we will study a generalization of the $f$-divergence functional for the class of 3-convex functions. It is an extension of the results obtained in \cite{prvi3conv}. Throughout this section, when mentioning the interval $[m,M]$, we assume that $[m,M]\subseteq \mathbb{R}_+$. For a 3-convex function $f\colon [m,M]\to \mathbb{R}$ we give the following definition of generalized $f$-divergence functional found in \cite{prvi3conv}:
\begin{equation}\label{fdivg}
\tilde{D}_f(\boldsymbol{p},\boldsymbol{q})=\sum_{i=1}^nq_if\left(\frac{p_i}{q_i}\right).
\end{equation}

We can utilize Theorem \ref{tm2} to get an Edmundson-Lah-Ribari\v c type inequality for the above defined generalized $f$-divergence functional.

\begin{tm}\label{tm2pr}
Let $[m,M]\subset \mathbb{R}$ be an interval such that $m\le 1\le M$. Let $f$ be a 3-convex function on the interval $I$ whose interior contains $[m,M]$ and differentiable on $\langle m,M\rangle$. Let $\boldsymbol{p}=(p_1,...,p_n)$ and $\boldsymbol{p}=(q_1,...,q_n)$ be probability distributions such that $p_i/q_i\in [m,M]$ for every $i=1,...,n$. Then we have
 \begin{align}
  &\left(1-m\right)\left[\dfrac{f(M)-f(m)}{M-m}-\dfrac{f_+'(m)}{2}\right]-\dfrac{1}{2}\sum_{i=1}^n(p_i-mq_i)f'\left(\dfrac{p_i}{q_i}\right) \nonumber\\
\le  &\dfrac{M-1}{M-m}f(m)+\dfrac{1-m}{M-m}f(M)-\tilde{D}_f(\boldsymbol{p},\boldsymbol{q})\label{tm2rez1pr}\\
  \le  & \dfrac{1}{2}\sum_{i=1}^n(Mq_i-p_i)f'\left(\frac{p_i}{q_i}\right)- (M-1)\left[\dfrac{f(M)-f(m)}{M-m}-\dfrac{f_-'(M)}{2}\right]. \nonumber
 \end{align}
\end{tm}

\emph{Proof:} Let $\boldsymbol{x}=(x_1,...,x_n)$ such that $x_i\in [m,M]$ for $i=1,...,n$. Let $\phi$ be a 3-convex function on the interval $I$ whose interior contains $[m,M]$ and differentiable on $\langle m,M\rangle$. In the relation (\ref{tm2rez1}) we can replace
$$f\longleftrightarrow \boldsymbol{x}, \ \ \mathrm{and} \ \  A(\boldsymbol{x})=\sum_{i=1}^np_ix_i.$$
In that way we get 
\begin{align*}
 &(\bar{x}-m)\left[\dfrac{\phi(M)-\phi(m)}{M-m}-\dfrac{\phi_+'(m)}{2}\right]-\dfrac{1}{2}\sum_{i=1}^np_i(x_i-m)\phi'(x_i) \nonumber\\
\le  &\dfrac{M-\bar{x}}{M-m}\phi (m)+\dfrac{\bar{x}-m}{M-m}\phi (M)-\sum_{i=1}^np_i\phi (x_i)\\
  \le  & \dfrac{1}{2}\sum_{i=1}^np_i(M-x_i)\phi'(x_i)- (M-\bar{x})\left[\dfrac{\phi(M)-\phi(m)}{M-m}-\dfrac{\phi_-'(M)}{2}\right] \nonumber 
\end{align*}
where $\bar{x}=\sum_{i=1}^np_ix_i$.
Since the function $f$ satisfies the same assumtions as $\phi$, in the previous relation we can set
$$\phi=f, \ \ p_i=q_i \ \ \mathrm{and} \ \ x_i=\dfrac{p_i}{q_i},$$
and after calculating 
$$\bar{x}=\sum_{i=1}^nq_i\dfrac{p_i}{q_i}=\sum_{i=1}^np_i=1$$
we get (\ref{tm2rez1pr}).
\begin{flushright}
 $\Box$
\end{flushright}

By utilizing Theorem \ref{tm3} in the analogous way as above, we get a different Edmundson-Lah-Ribari\v c type inequality for the generalized $f$-divergence functional (\ref{fdivg}), and it is given in the following theorem.

\begin{tm}\label{tm3pr}
Let $[m,M]\subset \mathbb{R}$ be an interval such that $m\le 1\le M$. Let $f$ be a 3-convex function on the interval $I$ whose interior contains $[m,M]$ and differentiable on $\langle m,M\rangle$. Let $\boldsymbol{p}=(p_1,...,p_n)$ and $\boldsymbol{p}=(q_1,...,q_n)$ be probability distributions such that $p_i/q_i\in [m,M]$ for every $i=1,...,n$. Then we have
 \begin{align}
       &(M-1)\left[f_-'(M)-\dfrac{f(M)-f(m)}{M-m}\right]-\dfrac{f_-''(M)}{2}\sum_{i=1}^n\dfrac{(Mq_i-p_i)^2}{q_i} \nonumber\\
\le  &\dfrac{M-1}{M-m}f (m)+\dfrac{1-m}{M-m}f(M)-\tilde{D}_f(\boldsymbol{p},\boldsymbol{q})\label{tm3rez1pr}\\
  \le  &(1-m)\left[\dfrac{f(M)-f(m)}{M-m}-f_+'(m)\right]-\dfrac{f_+''(m)}{2}\sum_{i=1}^n\dfrac{(p_i-mq_i)^2}{q_i}. \nonumber
 \end{align}
\end{tm}

\begin{rem}
Theorem \ref{tm1} can be in analogue way applied to generalized Csisz\' ar divergence functional, but since this application is already shown in \cite{prvi3conv}, we omit it. 
\end{rem}

\begin{primjer}\label{primjer}
Let $\boldsymbol{p}=(p_1,...,p_n)$ and $\boldsymbol{p}=(q_1,...,q_n)$ be probability distributions and let $[m,M]\subset \mathbb{R}$ be an interval such that $m\le 1\le M$ and $p_i/q_i\in [m,M]$ for every $i=1,...,n$.
\begin{itemize}
\item[$\triangleright$] \textbf{Kullback-Leibler divergence} of the probability distributions $\boldsymbol{p}$ and $\boldsymbol{q}$ is defined as
$$D_{KL}(\boldsymbol{p},\boldsymbol{q})=\sum_{i=1}^nq_i\log \dfrac{q_i}{p_i},$$
and the corresponding generating function is $f(t)=t\log t, t>0$. We can calculate $f'''(t)=-\frac{1}{t^2}<0$, so the function $-f(t)=-t\log t$ is 3-convex.
Now it is obvious that for the Kullback-Leibler divergence the inequalities (\ref{tm2rez1pr}) and (\ref{tm3rez1pr}) hold with reversed signs of inequality, with
$$f_+'(m)=\log m+1, \ \ f_-'(M)=\log M+1$$
and
$$f_+''(m)=\dfrac{1}{m}, \ \ f_-''(M)=\dfrac{1}{M}.$$

\item[$\triangleright$] \textbf{Hellinger divergence} of the probability distributions $\boldsymbol{p}$ and $\boldsymbol{q}$ is defined as
$$D_{H}(\boldsymbol{p},\boldsymbol{q})=\dfrac{1}{2}\sum_{i=1}^n(\sqrt{q_i}-\sqrt{p_i})^2,$$
and the corresponding generating function is $f(t)=\frac{1}{2}(1-\sqrt{t})^2, t>0$. We see that $f'''(t)=-\frac{3}{8}t^{-\frac{5}{2}}<0$, so the function $-f(t)=-\frac{1}{2}(1-\sqrt{t})^2$ is 3-convex. 
It is clear that for the Hellinger divergence the inequalities (\ref{tm2rez1pr}) and (\ref{tm3rez1pr}) hold with reversed signs of inequality, with
$$f_+'(m)=-\dfrac{1}{2\sqrt{m}}+\dfrac{1}{2}, \ \ f_-'(M)=-\dfrac{1}{2\sqrt{M}}+\dfrac{1}{2}$$
and
$$f_+''(m)=\dfrac{1}{4\sqrt{m^3}}, \ \ f_-''(M)=\dfrac{1}{4\sqrt{M^3}}.$$

\item[$\triangleright$] \textbf{Renyi divergence} of the probability distributions $\boldsymbol{p}$ and $\boldsymbol{q}$ is defined as
$$D_{\alpha}(\boldsymbol{p},\boldsymbol{q})=\sum_{i=1}^nq_i^{\alpha-1}p_i^{\alpha}, \ \alpha \in \mathbb{R},$$
and the corresponding generating function is $f(t)=t^{\alpha}, t>0$. We calculate that $f'''(t)=\alpha (\alpha-1)(\alpha-2)t^{\alpha-3}$ and see that the function $f(t)=t^{\alpha}$ is 3-convex for $0\le \alpha \le 1$ and $\alpha \ge 2$, and $-f(t)=-t^{\alpha}$ is 3-convex for $\alpha\le 0$ and $1<\alpha<2$, and we have
$$f_+'(m)=\alpha m^{\alpha-1}, \ \ f_-'(M)=\alpha M^{\alpha-1},$$ 
$$f_+''(m)=\alpha (\alpha-1)m^{\alpha-2} \ \ \mathrm{and} \ \ f_-''(M)=\alpha (\alpha-1)M^{\alpha-2}.$$
As regards, the Renyi divergence, the inequalities (\ref{tm2rez1pr}) and (\ref{tm3rez1pr}) hold for $0\le \alpha \le 1$ and $\alpha \ge 2$, and if $\alpha\le 0$ or $1<\alpha<2$ the signs of inequality are reversed.

\item[$\triangleright$] \textbf{Harmonic divergence} of the probability distributions $\boldsymbol{p}$ and $\boldsymbol{q}$ is defined as
$$D_{Ha}(\boldsymbol{p},\boldsymbol{q})=\sum_{i=1}^n\dfrac{2p_iq_i}{p_i+q_i},$$
and the corresponding generating function is $f(t)=\frac{2t}{1+t}$. We can calculate $f'''(t)=\frac{12}{(1+t)^4}>0$, so the function $f$ is 3-convex. Now it is obvious that for the harmonic divergence the inequalities (\ref{tm2rez1pr}) and (\ref{tm3rez1pr}) hold with
$$f_+'(m)=\frac{2}{(1+m)^2}, \ \ f_-'(M)=\frac{2}{(1+M)^2}$$
and
$$f_+''(m)=-\frac{4}{(1+m)^3}, \ \ f_-''(M)=-\frac{4}{(1+M)^3}.$$

\item[$\triangleright$] \textbf{Jeffreys divergence} of the probability distributions $\boldsymbol{p}$ and $\boldsymbol{q}$ is defined as
$$D_{J}(\boldsymbol{p},\boldsymbol{q})=\dfrac{1}{2}\sum_{i=1}^n(q_i-p_i)\log \dfrac{q_i}{p_i},$$
and the corresponding generating function is $f(t)=(1-t)\log \frac{1}{t}, t>0$. We see that $f'''(t)=-\frac{1}{t^2}-\frac{2}{t^3}<0$, so the function $-f(t)=(1-t)\log t$ is 3-convex. Instantly we get that for the Jeffreys divergence the inequalities (\ref{tm2rez1pr}) and (\ref{tm3rez1pr}) hold with reversed signs of inequality, with
$$f_+'(m)=\log m-\dfrac{1}{m}+1, \ \ f_-'(M)=\log M-\dfrac{1}{M}+1$$
and
$$f_+''(m)=\dfrac{1}{m}+\dfrac{1}{m^2}, \ \ f_-''(M)=\dfrac{1}{M}+\dfrac{1}{M^2}.$$
\end{itemize}
\end{primjer}

\section{Examples with Zipf and Zipf-Mandelbrot law}

Zipf’s law \cite{zipf}, \cite{zipf1} has and continues to attract considerable attention in a wide variety of scientific disciplines - from astronomy to demographics to software
structure to economics to zoology, and even to warfare \cite{war}. It is one of the basic laws in information science and bibliometrics, but it is also often used in linguistics. Same law in mathematical sense is also used in other scientific disciplines, but name of the law can be different, since regularities in different scientific fields are discovered independently from each other. Typically one is dealing with integer-valued observables (numbers of objects, people, cities, words, animals, corpses) and the frequency of their occurrence.

Probability mass function of Zipf's law with parameters $N\in \mathbb{N}$ and $s>0$ is:
\begin{equation*}
f(k;N,s)=\frac{1/k^s}{H_{N,s}}, \ \ \mathrm{where} \ \ H_{N,s}=\sum_{i=1}^N\frac{1}{i^s}.
\end{equation*}

Benoit Mandelbrot in 1966 gave an improvement of Zipf law for the count of the low-rank words. Various scientific fields use this law for different purposes, for example information sciences use it for indexing \cite{egg,sila}, ecological field studies in predictability of ecosystem \cite{moi}, in music it is used to determine aesthetically pleasing music \cite{mana}.

Zipf–Mandelbrot law is a discrete probability distribution
with parameters $N\in \mathbb{N}$, $q,s\in \mathbb{R}$ such that $q\ge 0$ and $s>0$, possible values $\{1, 2, ..., N\}$ and probability mass function
\begin{equation}\label{ZM}
f(i;N,q,s)=\frac{1/(i+q)^s}{H_{N,q,s}}, \ \ \mathrm{where} \ \ H_{N,q,s}=\sum_{i=1}^N\frac{1}{(i+q)^s}.
\end{equation}

Let $\boldsymbol{p}$ and $\boldsymbol{q}$ be Zipf-Mandelbrot laws with parameters $N\in \mathbb{N}$, $q_1,q_2\ge 0$ and $s_1,s_2>0$ respectively and let us denote
\begin{align}
m_{\boldsymbol{p},\boldsymbol{q}}&:=\mathrm{min}\left\{\dfrac{p_i}{q_i}\right\}=\dfrac{H_{N,q_2,s_2}}{H_{N,q_1,s_1}}\mathrm{min}\left\{\dfrac{(i+q_2)^{s_2}}{(i+q_1)^{s_1}}\right\}\nonumber\\
M_{\boldsymbol{p},\boldsymbol{q}}&:=\mathrm{max}\left\{\dfrac{p_i}{q_i}\right\}=\dfrac{H_{N,q_2,s_2}}{H_{N,q_1,s_1}}\mathrm{max}\left\{\dfrac{(i+q_2)^{s_2}}{(i+q_1)^{s_1}}\right\}\label{defmMKL}
\end{align}

In this section we utilize the results regarding Csisz\' ar divergence from the previous section in order to obtain different inequalities for the Zipf-Mandelbrot law. The first result that follows is a special case of Theorem \ref{tm2pr}, and it gives us Edmundson-Lah-Ribari\v c type inequality for the generalized $f$-divergence of the Zipf–Mandelbrot law.

\begin{kor}\label{kor2ZM}
Let $\boldsymbol{p}$ and $\boldsymbol{q}$ be Zipf-Mandelbrot laws with parameters $N\in \mathbb{N}$, $q_1,q_2\ge 0$ and $s_1,s_2>0$ respectively, and let $m_{\boldsymbol{p},\boldsymbol{q}}$ and $M_{\boldsymbol{p},\boldsymbol{q}}$ be defined in (\ref{defmMKL}). Let $f\colon [m_{\boldsymbol{p},\boldsymbol{q}},M_{\boldsymbol{p},\boldsymbol{q}}]\to \mathbb{R}$ be a 3-convex function. Then we have
 \begin{align}
  &\left(1-m_{\boldsymbol{p},\boldsymbol{q}}\right)\left[\dfrac{f(M_{\boldsymbol{p},\boldsymbol{q}})-f(m_{\boldsymbol{p},\boldsymbol{q}})}{M_{\boldsymbol{p},\boldsymbol{q}}-m_{\boldsymbol{p},\boldsymbol{q}}}-\dfrac{f_+'(m_{\boldsymbol{p},\boldsymbol{q}})}{2}\right]\nonumber\\
  & \ \ -\dfrac{1}{2}\sum_{i=1}^n\left(\dfrac{1}{(i+q_1)^{s_1}H_{N,q_1,s_1}}-\dfrac{m_{\boldsymbol{p},\boldsymbol{q}}}{(i+q_2)^{s_2}H_{N,q_2,s_2}}\right)f'\left(\dfrac{H_{N,q_2,s_2}}{H_{N,q_1,s_1}}\dfrac{(i+q_2)^{s_2}}{(i+q_1)^{s_1}}\right) \nonumber\\
\le  &\dfrac{M_{\boldsymbol{p},\boldsymbol{q}}-1}{M_{\boldsymbol{p},\boldsymbol{q}}-m_{\boldsymbol{p},\boldsymbol{q}}}f(m_{\boldsymbol{p},\boldsymbol{q}})+\dfrac{1-m_{\boldsymbol{p},\boldsymbol{q}}}{M_{\boldsymbol{p},\boldsymbol{q}}-m_{\boldsymbol{p},\boldsymbol{q}}}f(M_{\boldsymbol{p},\boldsymbol{q}})-\tilde{D}_f(\boldsymbol{p},\boldsymbol{q})\label{kor2rez1ZM}\\
  \le  & \dfrac{1}{2}\sum_{i=1}^n\left(\dfrac{M_{\boldsymbol{p},\boldsymbol{q}}}{(i+q_2)^{s_2}H_{N,q_2,s_2}}-\dfrac{1}{(i+q_1)^{s_1}H_{N,q_1,s_1}}\right)f'\left(\dfrac{H_{N,q_2,s_2}}{H_{N,q_1,s_1}}\dfrac{(i+q_2)^{s_2}}{(i+q_1)^{s_1}}\right)\nonumber\\
  & \ \ -(M_{\boldsymbol{p},\boldsymbol{q}}-1)\left[\dfrac{f(M_{\boldsymbol{p},\boldsymbol{q}})-f(m_{\boldsymbol{p},\boldsymbol{q}})}{M_{\boldsymbol{p},\boldsymbol{q}}-m_{\boldsymbol{p},\boldsymbol{q}}}-\dfrac{f_-'(M_{\boldsymbol{p},\boldsymbol{q}})}{2}\right]. \nonumber
 \end{align}
\end{kor}

Next result follows directly from Theorem \ref{tm3}, and it gives us another Edmundson-Lah-Ribari\v c type inequality for the generalized $f$-divergence of the Zipf–Mandelbrot law.

\begin{kor}\label{kor3ZM}
Let $\boldsymbol{p}$ and $\boldsymbol{q}$ be Zipf-Mandelbrot laws with parameters $N\in \mathbb{N}$, $q_1,q_2\ge 0$ and $s_1,s_2>0$ respectively, and let $m_{\boldsymbol{p},\boldsymbol{q}}$ and $M_{\boldsymbol{p},\boldsymbol{q}}$ be defined in (\ref{defmMKL}). Let $f\colon [m_{\boldsymbol{p},\boldsymbol{q}},M_{\boldsymbol{p},\boldsymbol{q}}]\to \mathbb{R}$ be a 3-convex function. Then we have
 \begin{align}
     &(M_{\boldsymbol{p},\boldsymbol{q}}-1)\left[f_-'(M_{\boldsymbol{p},\boldsymbol{q}})-
     \dfrac{f(M_{\boldsymbol{p},\boldsymbol{q}})-f(m_{\boldsymbol{p},\boldsymbol{q}})}{M_{\boldsymbol{p},\boldsymbol{q}}-m_{\boldsymbol{p},\boldsymbol{q}}}\right]\nonumber\\
     & \ \ -\dfrac{f_-''(M_{\boldsymbol{p},\boldsymbol{q}})}{2}\sum_{i=1}^n(i+q_2)^{s_2}H_{N,q_2,s_2}\left(\dfrac{M_{\boldsymbol{p},\boldsymbol{q}}}{(i+q_2)^{s_2}H_{N,q_2,s_2}}-\dfrac{1}{(i+q_1)^{s_1}H_{N,q_1,s_1}}\right)^2 \nonumber\\
\le  &\dfrac{M_{\boldsymbol{p},\boldsymbol{q}}-1}{M_{\boldsymbol{p},\boldsymbol{q}}-m_{\boldsymbol{p},\boldsymbol{q}}}f(m_{\boldsymbol{p},\boldsymbol{q}})+\dfrac{1-m_{\boldsymbol{p},\boldsymbol{q}}}{M_{\boldsymbol{p},\boldsymbol{q}}-m_{\boldsymbol{p},\boldsymbol{q}}}f(M_{\boldsymbol{p},\boldsymbol{q}})-\tilde{D}_f(\boldsymbol{p},\boldsymbol{q})\label{kor3rez1ZM}\\
  \le  &(1-m_{\boldsymbol{p},\boldsymbol{q}})\left[\dfrac{f(M_{\boldsymbol{p},\boldsymbol{q}})-f(m_{\boldsymbol{p},\boldsymbol{q}})}{M_{\boldsymbol{p},\boldsymbol{q}}-m_{\boldsymbol{p},\boldsymbol{q}}}-f_+'(m_{\boldsymbol{p},\boldsymbol{q}})\right]\nonumber\\
  & \ \ -\dfrac{f_+''(m_{\boldsymbol{p},\boldsymbol{q}})}{2}\sum_{i=1}^n(i+q_2)^{s_2}H_{N,q_2,s_2}\left(\dfrac{1}{(i+q_1)^{s_1}H_{N,q_1,s_1}}-\dfrac{m_{\boldsymbol{p},\boldsymbol{q}}}{(i+q_2)^{s_2}H_{N,q_2,s_2}}\right)^2. \nonumber
 \end{align}
\end{kor}

\begin{rem}
By taking into consideration Example \ref{primjer} one can see that Corollary \ref{kor2ZM} and Corollary \ref{kor3ZM} can easily be applied to any of the following divergences: Kullback-Leibler divergence, Hellinger divergence, Renyi divergence, harmonic divergence or Jeffreys divergence.
\end{rem}

\section{Exponential convexity}

Let $L$ satisfy conditions (L1) and (L2) on a non-empty set $E$ and let $A$ be any  positive linear functional on $L$ with $A(\boldsymbol{1})=1$. Let $\phi$ be a 3-convex function defined on an interval of real numbers $I$ whose interior contains the interval $[m,M]$ and let $f \in L$ such that $ \phi \circ f \in L$.

Motivated by inequalities (\ref{tm1rez1}), we define following linear functionals which represent the difference between the right and the left sides of the mentioned inequalities:
\begin{align}
    \Gamma_1(\phi)&=\dfrac{M-A(f)}{M-m}\phi (m)+\dfrac{A(f)-m}{M-m}\phi (M)-A\left(\phi (f)\right)\nonumber\\
    & \ \ \ -\dfrac{A\left[(M\boldsymbol{1}-f)(f-m\boldsymbol{1})\right]}{M-m}\left(\dfrac{\phi(M)-\phi(m)}{M-m}-\phi'(m)\right)\label{funkc1tm1}
    \end{align}
\begin{align}    
    \Gamma_2(\phi)&=\dfrac{A\left[(M\boldsymbol{1}-f)(f-m\boldsymbol{1})\right]}{M-m}\left(\phi'(M)-\dfrac{\phi(M)-\phi(m)}{M-m}\right)\nonumber\\
    & \ \ \ -\dfrac{M-A(f)}{M-m}\phi (m)+\dfrac{A(f)-m}{M-m}\phi (M)-A\left(\phi (f)\right)\label{funkc2tm1}
\end{align}

From Theorem \ref{tm1} it follows that functionals $\Gamma_1$ and $\Gamma_2$ are positive linear functionals under aforementioned assumptions.

If function $\phi$ is in addition differentiable on $\langle m,M\rangle$, then motivated by series of inequalities (\ref{tm2rez1}) and (\ref{tm3rez1}), we define following linear functionals:
\begin{align}
    \Gamma_3(\phi)&=\dfrac{M-A(f)}{M-m}\phi (m)+\dfrac{A(f)-m}{M-m}\phi (M)-A\left(\phi (f)\right)\nonumber\\
    & \ \ \ -(A(f)-m)\left[\dfrac{\phi(M)-\phi(m)}{M-m}-\dfrac{\phi'(m)}{2}\right]-\dfrac{1}{2}A[(f-m\boldsymbol{1})\phi'(f)] \label{funkc1tm2}
        \end{align}
\begin{align} 
    \Gamma_4(\phi)&=\dfrac{1}{2}A[(M\boldsymbol{1}-f)\phi'(f)]- (M-A(f))\left[\dfrac{\phi(M)-\phi(m)}{M-m}-\dfrac{\phi'(M)}{2}\right]\nonumber\\
    & \ \ \ -\dfrac{M-A(f)}{M-m}\phi (m)+\dfrac{A(f)-m}{M-m}\phi (M)-A\left(\phi (f)\right).\label{funkc2tm2}
\end{align}
\begin{align}
    \Gamma_5(\phi)&=\dfrac{M-A(f)}{M-m}\phi (m)+\dfrac{A(f)-m}{M-m}\phi (M)-A\left(\phi (f)\right)\nonumber\\
    & \ \ \ -(M-A(f))\left[\phi'(M)-\dfrac{\phi(M)-\phi(m)}{M-m}\right]-\dfrac{\phi''(M)}{2}A[(M\boldsymbol{1}-f)^2] \label{funkc1tm3}
        \end{align}
\begin{align} 
    \Gamma_6(\phi)&=(A(f)-m)\left[\dfrac{\phi(M)-\phi(m)}{M-m}-\phi'(m)\right]-\dfrac{\phi''(m)}{2}A[(f-m\boldsymbol{1})^2]\nonumber\\
    & \ \ \ -\dfrac{M-A(f)}{M-m}\phi (m)+\dfrac{A(f)-m}{M-m}\phi (M)-A\left(\phi (f)\right).\label{funkc2tm3}
\end{align}

Functionals $\Gamma_3$, $\Gamma_4$, $\Gamma_5$ and $\Gamma_6$ respectively represent the difference between the right and the left sides of the aforementioned inequalities and from Theorem \ref{tm2} and \ref{tm3} it follows that under the assumptions from above they are positive.

Now, let $[m,M]\subset \mathbb{R}$ be an interval such that $m\le 1\le M$. Let $f$ be a 3-convex function on the interval $I$ whose interior contains $[m,M]$ and differentiable on $\langle m,M\rangle$. Let $\boldsymbol{p}=(p_1,...,p_n)$ and $\boldsymbol{p}=(q_1,...,q_n)$ be probability distributions such that $p_i/q_i\in [m,M]$ for every $i=1,...,n$. Following linear functionals arise from series of inequalities (\ref{tm2rez1pr}) and (\ref{tm3rez1pr}) :
\begin{align}
       \Gamma_7(f)&=\dfrac{M-1}{M-m}f(m)+\dfrac{1-m}{M-m}f(M)-\tilde{D}_f(\boldsymbol{p},\boldsymbol{q})\nonumber\\
   & \ \ \ \ \ -\left(1-m\right)\left[\dfrac{f(M)-f(m)}{M-m}-\dfrac{f_+'(m)}{2}\right]-\dfrac{1}{2}\sum_{i=1}^n(p_i-mq_i)f'\left(\dfrac{p_i}{q_i}\right) \label{funkc1tm2pr}
\end{align}
\begin{align}
 \Gamma_8(f)&=\dfrac{1}{2}\sum_{i=1}^n(Mq_i-p_i)f'\left(\frac{p_i}{q_i}\right)- (M-1)\left[\dfrac{f(M)-f(m)}{M-m}-\dfrac{f_-'(M)}{2}\right]\nonumber\\
 & \ \ \ \ \ -\dfrac{M-1}{M-m}f(m)+\dfrac{1-m}{M-m}f(M)-\tilde{D}_f(\boldsymbol{p},\boldsymbol{q})\label{funkc2tm2pr}
\end{align}

\begin{align}
    \Gamma_9(f)&=\dfrac{M-1}{M-m}f(m)+\dfrac{1-m}{M-m}f(M)-\tilde{D}_f(\boldsymbol{p},\boldsymbol{q})\nonumber\\
   & \ \ \ \ \ - (M-1)\left[f_-'(M)-\dfrac{f(M)-f(m)}{M-m}\right]-\dfrac{f_-''(M)}{2}\sum_{i=1}^n\dfrac{(Mq_i-p_i)^2}{q_i}\label{funkc1tm3pr} 
\end{align}

\begin{align}
   \Gamma_{10}(f)&= (1-m)\left[\dfrac{f(M)-f(m)}{M-m}-f_+'(m)\right]-\dfrac{f_+''(m)}{2}\sum_{i=1}^n\dfrac{(p_i-mq_i)^2}{q_i}\nonumber\\
    & \ \ \ \ \ - \dfrac{M-1}{M-m}f (m)+\dfrac{1-m}{M-m}f(M)-\tilde{D}_f(\boldsymbol{p},\boldsymbol{q})\label{funkc2tm3pr} 
\end{align}
where $\tilde{D}_f(\boldsymbol{p},\boldsymbol{q})$ is the generalized f-divergence of the probability distributions $\boldsymbol{p}$ and $\boldsymbol{q}$ defined in (\ref{fdivg}). From Theorem \ref{tm2pr} and \ref{tm3pr} it immediately follows that $\Gamma_i$, $i=7,...,10$ are positive linear functionals on the class of 3-convex functions on the interval $I$ whose interior contains $[m,M]$ that are differentiable on $\langle m,M\rangle$.

First we will give some definitions and basic results regarding the exponential convexity that we need in the rest of this section. For the rest of this section $I$ will denote an interval of real numbers.

\begin{df}\label{nec}
A function $f\colon I\to \mathbb{R}$ is said to be $n$-exponentially convex in the Jensen sense on $I$ if
$$\sum\limits_{i,j=1}^n {\xi_i
\xi_j f \left(\frac{t_i + t_j}{2} \right)} \ge 0$$
holds for all choices of $\xi_i \in \mathbb{R}$ and every $t_i \in I$, $ i=1 ,... , n$.\\
A function $f\colon I\to \mathbb{R}$ is $n$-exponentially convex if it is $n$-exponentially convex in the Jensen sense and continuous on $I$.
\end{df}

\begin{rem}
It is clear from the definition that 1-exponentially convex functions in the Jensen sense are in fact nonnegative functions.
Also, $n$-exponentially convex functions in the Jensen sense are $k$-exponentially convex in the Jensen sense for every $k\in \mathbb{N}, \,\, k\le n$.
\end{rem}

\begin{df}\label{ec}
A function $f\colon I\to \mathbb{R}$ is exponentially convex in the Jensen sense on $I$ if
it is $n$-exponentially convex in the Jensen sense for all $n\in \mathbb{N}$.\\
A function $f\colon I\to \mathbb{R}$ is exponentially convex if it is exponentially convex in the Jensen sense and continuous on $I$.
\end{df}

\begin{rem}\label{remark1}
It is known that $f\colon I\to \mathbb{R}^+$ is  $\log$-convex in the Jensen sense, i.e.
\begin{align}
f\Big(\dfrac{t_1+t_2}{2}\Big)\le f(t_1)f(t_2) \ \ \mathrm{for \ all} \ \ t_1,t_2\in I\label{log}
\end{align}
 if and only if
 $$l^{2}f(t_1)+2lmf\left(\frac{t_1+t_2}{2}\right)+m^{2}f(t_2)\ge0$$
holds for each $l, m\in\mathbb{R}$ and $t_1, t_2 \in I$.
\end{rem}

Next results follow on the basis of the method developed by Jak\v seti\' c et. al. in \cite{julije}, where it is shown how positive linear functionals can be used to construct some new families of exponentially convex functions, so we give them with proofs omitted (see also \cite{ELRpaper}, \cite{roq}).

\begin{tm}\label{tmec}
Let $\Gamma_i$, $i=1,...,10$ be linear functionals defined in (\ref{funkc1tm1})-(\ref{funkc2tm3pr}) respectively, with corresponding assumptions.
Let $J$ be an interval in $\mathbb{R}$, and let $\Upsilon=\{\phi_t\colon [m,M]\to \mathbb{R} | t\in J \}$ be a family of differentiable functions such that for every four distinct points $u_0,u_1,u_2,u_3\in [m,M]$ the mapping $t\longmapsto [u_0,u_1,u_2,u_3]\phi_t$ is $n$-exponentially convex in the Jensen sense. Then the mapping $t\longmapsto \Gamma_i (\phi_t)$ is $n$-exponentially convex in the Jensen sense on $J$ for $i=1,...,10$. If additionally the mapping $t\longmapsto \Gamma_i (\phi_t)$ is continuous on $J$ for $i=1,...,10$, then it is $n$-exponentially convex on $J$. 
\end{tm}

If the assumptions of Theorem \ref{tmec} hold for all $n\in \mathbb{N}$, then we immediately get the following corollary.

\begin{kor}\label{corec}
Let $\Gamma_i$, $i=1,...,10$ be linear functionals defined in (\ref{funkc1tm1})-(\ref{funkc2tm3pr}) respectively, with corresponding assumptions.
Let $J$ be an interval in $\mathbb{R}$, and let $\Upsilon=\{\phi_t\colon [m,M]\to \mathbb{R} | t\in J \}$ be a family of differentiable functions such that for every four distinct points $u_0,u_1,u_2,u_3\in [m,M]$ the mapping $t\longmapsto [u_0,u_1,u_2,u_3]\phi_t$ is exponentially convex in the Jensen sense. Then the mapping $t\longmapsto \Gamma_i (\phi_t)$ is exponentially convex in the Jensen sense on $J$ for $i=1,...,10$. If additionally the mapping $t\longmapsto \Gamma_i (\phi_t)$ is continuous on $J$ for $i=1,...,10$, then it is exponentially convex on $J$. 
\end{kor}

\begin{kor}\label{corec2}
Let $\Gamma_i$, $i=1,...,10$ be linear functionals defined in (\ref{funkc1tm1})-(\ref{funkc2tm3pr}) respectively, with corresponding assumptions.
Let $J$ be an interval in $\mathbb{R}$, and let $\Upsilon=\{\phi_t\colon [m,M]\to \mathbb{R} | t\in J \}$ be a family of differentiable functions such that for every four distinct points $u_0,u_1,u_2,u_3\in [m,M]$ the mapping $t\longmapsto [u_0,u_1,u_2,u_3]\phi_t$ is $2$-exponentially convex in the Jensen sense. Then the following statements hold.
\begin{itemize}
\item[(i)] If the mapping $t\longmapsto \Gamma_i (\phi_t)$ is continuous on $J$, then for $r,s,t\in J$ such that $r<s<t$ we have
\begin{align}
\Gamma_i (\phi_s)^{t-r}\le \Gamma_i (\phi_r)^{t-s}\Gamma_i (\phi_t)^{s-r}\label{corec2a}
\end{align}
for $i=1,...,10$.
\item[(ii)] If the mapping $t\longmapsto \Gamma_i (\phi_t)$ is strictly positive and differentiable on $J$, then for all $s,t,u,v\in J$ such that $s\le u$ and $t\le v$ we have
\begin{align*}
\mathfrak{B}_{s,t}(\Upsilon)\le \mathfrak{B}_{u,v}(\Upsilon),
\end{align*}
where
\begin{equation}\label{corec2b}
\mathfrak{B}_{s,t}(\Upsilon)= \left\{
           \begin{array}{ll}
            \left(\dfrac{\Gamma_i(\phi_s)}{\Gamma_i(\phi_t)}\right)^{\frac{1}{s-t}}, &s\neq t\\
             \exp \left(\dfrac{\frac{d}{d{s}}(\Gamma_i(\phi_s))}{\Gamma_i(\phi_s)}\right),  & s=t
           \end{array}
         \right.
\end{equation}
for $i=1,...,10$.
\end{itemize}
\end{kor}

\section{Stolarsky-type means}

First we will give two mean value results, which are essential in producing criteria under which Stolarsky-type quotients are actual means. The results below are proven by following the steps in the proof of corresponding theorems from \cite{julije}, so we omit the proof.

\begin{tm}\label{tm_srednja1}
Let $\Gamma_i$, $i=1,...,10$ be linear functionals defined in (\ref{funkc1tm1})-(\ref{funkc2tm3pr}) respectively with corresponding assumptions. Then for $\phi \in \mathcal{C}^3([m,M])$ there exists $\xi \in [m,M]$ such that
\begin{align*}
\Gamma_i (\phi)&=\dfrac{\phi'''(\xi)}{6}\Gamma_i (\phi_0)
\end{align*}
for $i=1,...,10$, where $\phi_0(t)=t^3$.
\end{tm}

\begin{tm}\label{tm_srednja2}
Let $\Gamma_i$, $i=1,...,10$ be linear functionals defined in (\ref{funkc1tm1})-(\ref{funkc2tm3pr}) respectively with corresponding assumptions. Let $\phi_1,\phi_2 \in \mathcal{C}^3([m,M])$. If $\Gamma_i (\phi_2)\neq 0$, then there exists $\xi \in [m,M]$ such that 
$$\dfrac{\Gamma_i (\phi_1)}{\Gamma_i (\phi_2)}=\dfrac{\phi_1'''(\xi)}{\phi_2'''(\xi)} \ \ \mathrm{for} \ \ i=1,...,10$$
or
$$\phi_1'''(\xi)=\phi_2'''(\xi)=0.$$
\end{tm}

\begin{rem}
If the inverse of the function $\dfrac{\phi_1'''}{\phi_2'''}$ exists, then various kinds of means can be defined by Theorem \ref{tm_srednja2}. That is,
\begin{equation}\label{mean}
    \xi=\left(\dfrac{\phi_1'''}{\phi_2'''}\right)^{-1}\left(\dfrac{\Gamma_i (\phi_1)}{\Gamma_i (\phi_2)}\right)
\end{equation}
for $i=1,...,10$.
\end{rem}

Let us consider the following family of functions
$$\Upsilon_1=\{\phi_t\colon [m,M]\to \mathbb{R} \ | \ t\in \mathbb{R}\}, \ \ \ [m,M]\subset \langle 0,+\infty \rangle,$$
defined by
\begin{equation}\label{funkcije1}
\phi_t(x)= \left\{
           \begin{array}{llll}
            \dfrac{1}{t(t-1)(t-2)}x^t, &t\neq 0,1,2\\
             \dfrac{1}{2}\ln x,  & t=0,\\
             -x\ln x, &t=1,\\
            \dfrac{1}{2}x^2\ln x, &t=2.
           \end{array}
         \right.
\end{equation}

Since $\phi_t'''(x)=x^{t-3}\ge 0$, the functions $\phi_t$ are $3$-convex, and the function
$$\phi(x)=\sum_{i,j=1}^n\xi_i\xi_j\phi_{\frac{t_i+t_j}{2}}(x)$$
satisfies
$$\phi'''(x)=\sum_{i,j=1}^n\xi_i\xi_j\phi_{\frac{t_i+t_j}{2}}'''(x)=\Big(\sum_{i=1}^n\xi_ie^{\left(\frac{t_i}{2}-3\right)\ln x}\Big)^2\ge 0,$$
so $\phi$ is 3-convex. Therefore we have
$$0\le [u_0,u_1,u_2,u_3]\phi=\sum_{i,j=1}^n\xi_i\xi_j[u_0,u_1,u_2,u_3]\phi_{\frac{t_i+t_j}{2}}(x),$$
so the mapping $t\longmapsto [u_0,u_1,u_2,u_3]\phi_t$ is $n$-exponentially convex in the Jensen sense. Since this holds for every $n\in \mathbb{N}$, we see that family $\Upsilon_1$ satisfies the assumptions of Corollary \ref{corec}. Hence, the mapping $t\longmapsto \Gamma_i (\phi_t)$ is exponentially convex in the Jensen sense. It is easy to check that that it is also continuous, so the mappings $t\longmapsto \Gamma_i (\phi_t)$, $i=1,...,10$, are exponentially convex.

If we apply Theorem \ref{tm_srednja2} for functions $\phi_1=\phi_s$ and $\phi_2=\phi_t$ given by (\ref{funkcije1}), we can conclude that there exists $\xi \in [m,M]\subset \langle 0,+\infty \rangle$ such that
$$\xi=\Big(\dfrac{\phi_s'''}{\phi_t'''}\Big)^{-1}\Big(\dfrac{\Gamma_i (\phi_s)}{\Gamma_i (\phi_t)}\Big)=\Big(\dfrac{\Gamma_i (\phi_s)}{\Gamma_i (\phi_t)}\Big)^{\frac{1}{s-t}}, \ \ s\neq t.$$

Therefore, $\mathfrak{B}_{s,t}(\Upsilon_1)$ given by (\ref{corec2b}) for the family of functions $\Upsilon_1$ is a mean of the segment $[m,M]$. The limiting cases $s\rightarrow t$ can be calculated, and are equal to:
\begin{equation*}
\mathfrak{B}_{s,t}(\Upsilon_1)= \left\{
           \begin{array}{lllll}
            \left(\dfrac{\Gamma_i(\phi_s)}{\Gamma_i(\phi_t)}\right)^{\frac{1}{s-t}}, &s\neq t,\\
             \exp \left(\dfrac{2\Gamma_i(\phi_s\phi_0)}{\Gamma_i(\phi_0)}-\dfrac{3s^2-6s+2}{s(s-1)(s-2)}\right),  & s=t\neq 0,1,2,\\
              \exp \left(\dfrac{\Gamma_i(\phi^2_0)}{\Gamma_i(\phi_0)}+\dfrac{3}{2}\right),  & s=t=0,\\
              \exp \left(\dfrac{\Gamma_i(\phi_0\phi_1)}{\Gamma_i(\phi_1)}\right),  & s=t=1,\\
              \exp \left(\dfrac{\Gamma_i(\phi_0\phi_2)}{\Gamma_i(\phi_2)}-\dfrac{3}{2}\right),  & s=t=2.
           \end{array}
         \right.
\end{equation*}
for $i=1,...,10$. From Corollary \ref{corec2}(ii) it follows that the means $\mathfrak{B}_{s,t}(\Upsilon_1)$ are monotone in parameters $s$ and $t$.

Now consider a family of functions
$$\Upsilon_2=\{\varphi_t\colon [m,M]\to \mathbb{R} \ | \ t\in \mathbb{R}\}, \ \ \ [m,M]\subset \langle 0,+\infty \rangle,$$
defined by
\begin{equation*}
\varphi_t(x)= \left\{
           \begin{array}{ll}
           \dfrac{1}{t^3}e^{tx}, &t\neq 0\\
             \dfrac{1}{6}x^3,  & t=0
           \end{array}
         \right.
\end{equation*}

By straightforward calculation we see that $\varphi_t'''(x)=e^{tx}\ge 0$, so it follows that the functions $\varphi_t$ are $3$-convex. The function defined by
$$\varphi(x)=\sum_{i,j=1}^n\xi_i\xi_j\varphi_{\frac{t_i+t_j}{2}}(x)$$
satisfies
$$\varphi'''(x)=\sum_{i,j=1}^n\xi_i\xi_j\varphi_{\frac{t_i+t_j}{2}}'''(x)=\Big(\sum_{i=1}^n\xi_ie^{\frac{t_i}{2}x}\Big)^2\ge 0,$$
so it is also 3-convex. Consequently it holds
$$0\le [u_0,u_1,u_2,u_3]\varphi=\sum_{i,j=1}^n\xi_i\xi_j[u_0,u_1,u_2,u_3]\varphi_{\frac{t_i+t_j}{2}}(x),$$
so the mapping $t\longmapsto [u_0,u_1,u_2,u_3]\varphi_t$ is $n$-exponentially convex in the Jensen sense. Because this holds for every $n\in \mathbb{N}$, the family $\Upsilon_2$ satisfies the assumptions of Corollary \ref{corec}. Therefore, the mapping $t\longmapsto \Gamma_i (\varphi_t)$ is exponentially convex in the Jensen sense. Since it is also continuous, the mappings $t\longmapsto \Gamma_i (\varphi_t)$, $i=1,...,10$, are exponentially convex.

If we put $\phi_1=\varphi_s$ and $\phi_2=\varphi_t$ in Theorem \ref{tm_srednja2}, we see that there has to exist $\xi \in [m,M]\subset \langle 0,+\infty \rangle$ such that
$$\xi=\Big(\dfrac{\varphi_s'''}{\varphi_t'''}\Big)^{-1}\Big(\dfrac{\Gamma_i (\varphi_s)}{\Gamma_i (\varphi_t)}\Big)=\frac{1}{s-t}\mathrm{ln}\left(\dfrac{\Gamma_i (\varphi_s)}{\Gamma_i (\varphi_t)}\right), \ \ s\neq t.$$

Consequently, $\mathfrak{M}_{s,t}(\Upsilon_2)$ defined by 
\begin{equation*}
\mathfrak{M}_{s,t}(\Upsilon_2)= \left\{
           \begin{array}{ll}
            \dfrac{1}{s-t}\mathrm{ln}\left(\dfrac{\Gamma_i(\phi_s)}{\Gamma_i(\phi_t)}\right), &s\neq t\\
             \dfrac{\frac{d}{d{s}}(\Gamma_i(\phi_s))}{\Gamma_i(\phi_s)},  & s=t
           \end{array}
         \right.
\end{equation*}
for the family of functions $\Upsilon_2$ is a mean of the segment $[m,M]$. The limiting cases $s\rightarrow t$ can be calculated, and are equal to:

\begin{equation*}
\mathfrak{M}_{s,t}(\Upsilon_2)= \left\{
           \begin{array}{lll}
             \dfrac{1}{s-t}\mathrm{ln}\left(\dfrac{\Gamma_i(\phi_s)}{\Gamma_i(\phi_t)}\right), &s\neq t\\
            \dfrac{\Gamma_i(\mathrm{id}\cdot \phi_s)}{\Gamma_i(\phi_s)}-\dfrac{3}{s}, &s=t\neq 0,\\
             \dfrac{\Gamma_i(\mathrm{id}\cdot \phi_0)}{4\Gamma_i(\phi_0)},  & s=t=0,
           \end{array}
         \right.
\end{equation*}
for $i=1,...,10$. Notice that this is a monotonic mean (in respect to parameters $s$ and $t$).

\end{document}